\documentclass[runningheads]{llncs}
\usepackage{amsfonts, amssymb, amsmath}
\usepackage{dsfont} 
\usepackage{nicefrac}
\usepackage{graphicx}
\usepackage{subcaption}

\newcommand{\Symm}[1]{\mathcal{S}^{#1}}

\DeclareMathOperator\conv{conv}

\begin{document}

\title{Numerical analysis of the convex relaxation of the barrier parameter functional of self-concordant barriers}

\author{Vitali Pirau\inst{1} \and
Roland	Hildebrand\inst{1,2}}
\authorrunning{Vitali Pirau and Roland Hildebrand}
%
\institute{Moscow Institute of Physics and Technology (National Research University) (MIPT), Institutskiy per., 9, Dolgoprudny, Moscow region, 141701 \and
Skolkovo Institute of Science and Technology, Skolkovo Innovation Center, Bolshoy Boulevard, 30, p.1, Moscow 121205, Russia}
\maketitle    

\begin{abstract}
Self-concordant barriers are essential for interior-point algorithms in conic programming. To speed up the convergence it is of interest to find a barrier with the lowest possible parameter for a given cone. The barrier parameter is a non-convex function on the set of self-concordant barriers on a given cone, and finding an optimal barrier amounts to solving a non-convex infinite-dimensional optimization problem. In this work we study the degradation of the optimal value of the problem when the problem is convexified, and provide an estimate of the accuracy of the convex relaxation. The amount of degradation can be computed by comparing a 1-parameter family of non-convex bodies in $\mathbb{R}^3$ with their convex hulls. Our study provides insight into the degree of non-convexity of the problem and opens up the possibility of constructing suboptimal barriers by solving the convex relaxation.
\end{abstract}

\section{Introduction}

Conic programming deals with the problem of minimizing a linear functional over the intersection of an affine subspace with a regular convex cone, i.e., a closed convex cone with non-empty interior and containing no lines. The standard methods to solve such a problem are interior-point methods, which rely on the availability of a logarithmically homogeneous self-concordant barrier, a function defined on the interior of the cone and having the following properties \cite{NesNem94}.

{\definition \label{def:lhsc_barrier} Let $K \subset \mathbb R^n$ be a regular convex cone. A \emph{self-concordant logarithmically homogeneous barrier} on $K$ with parameter $\nu$ is a $C^3$ function $F: K^o \to \mathbb R$ satisfying
\begin{itemize}
\item $F(\alpha x) = -\nu\log\alpha + F(x)$ for all $\alpha > 0$, $x \in K^o$ (logarithmic homogeneity),
\item $F''(x) \succ 0$ for all $x \in K^o$ (locally strong convexity),
\item $\lim_{x \to\partial K} F(x) = +\infty$ (barrier property),
\item $|F'''(x)[u,u,u]| \leq 2(F''(x)[u,u])^{3/2}$ for all $x \in K^o$, $u \in T_xK^o$ (self-concordance).
\end{itemize} }

The convergence speed of an interior-point method is determined by the value of the scalar parameter $\nu$ and increases with decreasing $\nu$. Well-known sub-classes of conic programs are linear programs (LP), second order cone programs (SOCP), and semi-definite programs (SDP), which are conic programs over symmetric cones. These optimization problems can be very efficiently solved due to the availability of computable self-concordant barriers with optimal (i.e., the lowest possible) barrier parameter \cite{NesNem94,nesterov1997self}. However, convex optimization problems whose constraints contain exponentials or polynomials other than linear or quadratic functions, e.g., geometric programs, lead to conic programs over non-symmetric cones. For most of the non-symmetric cones, among them the 3-dimensional power cones or the cones over the unit ball of the $p$-norms, barriers with an optimal parameter and even the optimal value of the parameter are not known. Instead, barriers are used which are constructed from the analytic description of the cone in question and are far from optimal \cite{nesterov2012towards}. Optimal barrier parameters are so far known only for homogeneous cones \cite{Guler98} and in the case when singularities on the cone boundary bound the parameter from below by the dimension of the cone \cite[Section 2.3.4]{NesNem94}.

In this contribution we advocate a new paradigm for non-symmetric cone programming. We consider the problem of finding an optimal barrier on a given cone as an optimization problem. This optimization problem is both non-convex and infinite-dimensional. While the second draw-back can be countered by a sequence of asymptotically exact finite-dimensional approximations of increasing complexity, the non-convexity is a more fundamental problem. However, we shall show that the convex relaxation of the problem leads to a controllable loss in terms of the objective value, i.e., the barrier parameter. More precisely, if $\nu_{rel}$ is the optimal value of the relaxed (convexified) problem, then the optimal value of the barrier parameter can be bounded by $\nu_{rel} \leq \nu_{opt} \leq \tilde\nu(\nu_{rel})$, where $\tilde\nu(\nu)$ is a piece-wise algebraic function. In the regime of small $\nu_{opt} \lesssim 4$, which is the relevant one for the low-dimensional cones arising from $p$-norm constraints and in power cone programming, the gap is also small and the relaxation becomes asymptotically exact if $\nu_{opt}$ tends to its lowest possible value 2.

The remainder of the paper is structured as follows. In Section \ref{sec:affine_sec} we use the property of logarithmic homogeneity in order to replace the interior of the cone as the domain of definition of the barrier by its intersection with an affine hyperplane. In Section \ref{sec:non_convex} we formulate the problem of minimizing the barrier parameter on a given cone as a non-convex infinite-dimensional optimization problem. In Section \ref{sec:numerical} we numerically examine the quality of the convex relaxation of this problem and derive a tight performance guaranty, calculating $\tilde\nu(\nu)$.

\section{Reduction to an affine section} \label{sec:affine_sec}

For a given cone $K$, a barrier $F$ on $K$ is defined on its interior $K^o$. However, on every ray in $K^o$ the value of $F$ at one point determines the values on the whole ray due to the logarithmic homogeneity of $F$. Therefore the values of $F$ on the interior of a compact section $C$ of $K$ by an affine hyperplane determine $F$ completely if the parameter $\nu$ is fixed. Logarithmic homogeneity can hence be used to lower the dimension of the domain of definition of $F$ by one and thus to simplify the considered optimization problem. We have the following result \cite[Theorem 2]{hildebrand2022barriers}.

{\lemma \label{lem:section} Let $K \subset R^{n+1}$ be a regular convex cone, $n \geq 1$, let $C$ be an $n$-dimensional compact proper affine section of $K$, and let $\nu$ be a positive number. To any $C^3$ function $f: C^o \to \mathbb R$ we associate a $C^3$ function $F: K^o \to \mathbb R$ by $F(\alpha x) = \nu(-\log\alpha + f(x))$ for all $x \in C^o$ and $\alpha > 0$. Then $F$ is a self-concordant logarithmically homogeneous barrier with parameter $\nu$ on $K$ if and only if $\nu \geq 2$ and $f$ satisfies the properties
\begin{itemize}
\item $f''(x) - f'(x) \otimes f'(x) \succ 0$ for all $x \in C^o$,
\item $\lim_{x \to \partial C} f(x) = +\infty$,
\item $|f'''(x)[v,v,v] - 6f''(x)[v,v]f'(x)[v] + 4(f'(x)[v])^3| \leq 2\gamma(f''(x)[v,v] - (f'(x)[v])^2)^{3/2}$ for all $x \in C^o$, $v \in T_xC^o$,
\end{itemize}
where $\gamma = \frac{\nu-2}{\sqrt{\nu-1}}$. }

Denote by $M_{\nu}$ the set of functions $f$ satisfying the conditions in Lemma \ref{lem:section}. This set is not convex for $\nu > 2$. We will show, however, that the convex hull of $M_{\nu}$ is contained in the set $M_{\nu'}$, where $\nu' > \nu$ is some number which depends only on $\nu$.

Note that the function $\gamma = \frac{\nu-2}{\sqrt{\nu-1}}$ is strictly monotonously increasing in $\nu$. Therefore the problem of minimizing $\nu$ such that a logarithmically homogeneous self-concordant barrier $F$ on $K$ with parameter $\nu$ exists is equivalent to the problem of minimizing $\gamma$ such that a $C^3$ function $f$ on $C^o$ satisfying the three conditions in Lemma \ref{lem:section} exists. In the sequel we shall consider exclusively the latter problem on compact convex sets $C$. However, we shall replace the last condition in Lemma \ref{lem:section} by a Lipschitz condition on the second derivative $f''$. This corresponds to passing to the closure of the set of barriers with a given parameter value $\nu$. This modification has no impact on the practical usability of the barrier, since interior-point methods require only the first two derivatives of the barrier. We thus look for functions satisfying the following definition.

{\definition \label{def:lipschitz} Let $C \subset R^n$ be a compact convex set, $n \geq 1$, let $\nu \geq 2$, and set $\gamma = \frac{\nu-2}{\sqrt{\nu-1}}$. We call a $C^2$ function $f: C^o \to \mathbb R$ \emph{admissible} if it satisfies the properties
\begin{itemize}
\item $f''(x) - f'(x) \otimes f'(x) \succ 0$ for all $x \in C^o$,
\item $\lim_{x \to \partial C} f(x) = +\infty$,
\item for all $x \in C^o$, $u \in T_xC^o$ we have
\end{itemize}
\begin{align*}
\limsup_{\epsilon \to 0} \frac{f''(x+\epsilon u)[u,u] - f''(x)[u,u]}{\epsilon} &\leq \zeta_+(x,u), \\
\liminf_{\epsilon \to 0} \frac{f''(x+\epsilon u)[u,u] - f''(x)[u,u]}{\epsilon} &\geq \zeta_-(x,u),
\end{align*}
where
\[ \zeta_{\pm}(x,u) = 6f''(x)[u,u]f'(x)[u] - 4(f'(x)[u])^3 \pm 2\gamma(f''(x)[u,u] - (f'(x)[u])^2)^{3/2}.
\]
}

The Lipschitz condition on $f''$ translates into a similar condition on the second derivative $F''$ of the barrier defined as in Lemma \ref{lem:section} from the function $f$. This Lipschitz condition on $F''$ replaces the stronger self-concordance condition in Definition \ref{def:lhsc_barrier}. 

\section{Reformulation as optimization problem} \label{sec:non_convex}

The conditions on the function $f$ in Definition \ref{def:lipschitz} are not well suited for optimization. In this section we shall reformulate these conditions in the form of constraints that are more common in optimization problems. To this end we shall first consider the case when the compact set $C$ is an interval.

\subsection{Admissible functions on an interval}

In this section we reformulate the conditions on a function $f: (-1,1) \to \mathbb R$ listed in Definition \ref{def:lipschitz}. Since these conditions are invariant under changes of $f$ by additive constants, we may ignore the values of $f$ and consider only the scalar functions $f',f''$ on $(-1,1)$. Let us first determine conditions on the pair $(f'(x),f''(x))$ at a single point $x \in (-1,1)$.

{\lemma \label{lem:dim1normalized} Let $f: (-1,1) \to \mathbb R$ be a $C^2$ function satisfying the properties in Definition \ref{def:lipschitz}. Then the pair $(f'(x),f''(x))$ satisfies the constraints
\[ \frac{\sqrt{f''(x)-f'(x)^2}}{\sqrt{\nu-1}} \leq f'(x) + \frac{1}{1+x} \leq \sqrt{\nu-1}\sqrt{f''(x)-f'(x)^2},
\]
\[ \frac{\sqrt{f''(x)-f'(x)^2}}{\sqrt{\nu-1}} \leq - f'(x) + \frac{1}{1-x} \leq \sqrt{\nu-1}\sqrt{f''(x)-f'(x)^2}
\]
for all $x \in (-1,1)$.}

\begin{proof}
    The second chain of inequalities follows from \cite[Corollary 6]{hildebrand2022barriers}. The first chain of inequalities follows from the second one by the change of variables $x \mapsto -x$.
\end{proof}

Further we have the following result \cite[Lemma 1]{hildebrand2022barriers}.

{\lemma \label{lem:2points} Let $f: (-1,1) \to \mathbb R$ be a $C^2$ function satisfying the properties in Definition \ref{def:lipschitz}. Let $x_0 \in (-1,1)$ and set $p_0 = f'(x_0)$, $h_0 = f''(x_0)$. Let $p_{\pm}(x)$ be the solutions of the differential equations
\begin{equation} \label{p_extr_diffeq}
p''_{\pm} = 6p'_{\pm}p_{\pm} - 4p_{\pm}^3 \pm 2\gamma(p'_{\pm} - p_{\pm}^2)^{3/2},
\end{equation}
respectively, with initial conditions $p_{\pm}(x_0) = p_0$, $p'_{\pm}(x_0) = h_0$. Then for every $x \in (-1,1)$ we have 
\[ p_-(x) \leq f'(x) \leq p_+(x),
\]
whenever the concerned functions are defined. }

Let us provide explicit expressions for the functions $p_{\pm}$. We have \cite[Lemma 1]{hildebrand2022barriers}
\begin{align*}
p_-(x) &= \frac{p_0 + (x-x_0)(g_0^2 - p_0^2 + \gamma g_0p_0)}{- g_0^2(x-x_0)^2 + (p_0(x-x_0) - 1)^2 - \gamma g_0(x-x_0)(p_0(x-x_0) - 1)}, \\
p_+(x) &= \frac{p_0 + (x-x_0)(g_0^2 - p_0^2 - \gamma g_0p_0)}{- g_0^2(x-x_0)^2 + (p_0(x-x_0) - 1)^2 + \gamma g_0(x-x_0)(p_0(x-x_0) - 1)},
\end{align*}
where $p_0 = p_{\pm}(x_0)$, $g_0 = \sqrt{h_0 - p_0^2}$.

We now show that the conditions in Lemmas \ref{lem:dim1normalized} and \ref{lem:2points} are not only necessary, but also sufficient to guarantee the admissibility of $f$ in the sense of Definition \ref{def:lipschitz}.

{\lemma \label{lem:sufficiency2dim} Let $\nu \geq 2$ and $\gamma = \frac{\nu-2}{\sqrt{\nu-1}}$, and let $\epsilon > 0$ be arbitrary. Let $p,h: (-1,1) \to \mathbb R$ be functions satisfying the properties
\begin{equation} \label{lemma4_1}
\frac{\sqrt{h(x)-p(x)^2}}{\sqrt{\nu-1}} \leq p(x) + \frac{1}{1+x} \leq \sqrt{\nu-1}\sqrt{h(x)-p(x)^2},    
\end{equation} 
\begin{equation} \label{lemma4_2}
\frac{\sqrt{h(x)-p(x)^2}}{\sqrt{\nu-1}} \leq - p(x) + \frac{1}{1-x} \leq \sqrt{\nu-1}\sqrt{h(x)-p(x)^2}
\end{equation}
for all $x \in (-1,1)$. Suppose further that for all $x_0 \in (-1,1)$ we have $p_-(x) \leq p(x) \leq p_+(x)$ for every $x$ which has distance less than $\epsilon$ to $x_0$ in the Hilbert metric of $(-1,1)$, where $p_{\pm}(x)$ are the solutions of \eqref{p_extr_diffeq} with initial conditions $p_{\pm}(x_0) = p(x_0)$, $p'_{\pm}(x_0) = h(x_0)$.

Then $p \in C^1$, $h \in C^0$, $p' = h$, and for every $x_0 \in (-1,1)$ we have
\begin{align*}
\limsup_{x \to x_0} \frac{h(x) - h(x_0)}{x-x_0} &\leq 6h(x_0)p(x_0) - 4p(x_0)^3 + 2\gamma(h(x_0) - p(x_0)^2)^{3/2}, \\
\liminf_{x \to x_0} \frac{h(x) - h(x_0)}{x-x_0} &\geq 6h(x_0)p(x_0) - 4p(x_0)^3 - 2\gamma(h(x_0) - p(x_0)^2)^{3/2}.
\end{align*}
}

\begin{proof}
First note that the upper and lower bounds $p_{\pm}$ of $p$ have the same value $p(x_0)$ and the same derivative $h(x_0)$ at $x_0$. Therefore $p(x)$ is continuous and differentiable at $x_0$ with $p'(x_0) = h(x_0)$.

Fix $x^* \in (-1,1)$. Since the functions $p$ and $h$ are locally bounded by virtue of the assumed inequalities there exists a constant $C > 0$ and a neighbourhood $U \subset (-1,1)$ of $x^*$ with diameter not exceeding $\epsilon$ such that $p_{\pm}$ are defined on $U$ for all initial points $x_0 \in U$ and the second derivatives $|p_{\pm}''|$ are uniformly bounded by $C$ on $U$ for all initial points $x_0 \in U$.

For $x_1 > x_0$, $x_0,x_1 \in U$ we get
\begin{align*}
    p(x_0) + h(x_0)(x_1-x_0) - \frac{C}{2}(x_1-x_0)^2 &\leq p_-(x_1) \leq p(x_1) \leq p_+(x_1) \\ &\leq p(x_0) + h(x_0)(x_1-x_0) + \frac{C}{2}(x_1-x_0)^2.
\end{align*} 
Exchanging the roles of $x_0,x_1$ we likewise obtain
\[ p(x_1) + h(x_1)(x_0-x_1) - \frac{C}{2}(x_0-x_1)^2 \leq p(x_0) \leq p(x_1) + h(x_1)(x_0-x_1) + \frac{C}{2}(x_0-x_1)^2.
\]
Equivalently we get
\[
\begin{aligned}
h(x_0) - \frac{C}{2}(x_1-x_0) &\leq& \frac{p(x_1)-p(x_0)}{x_1-x_0} &\leq& h(x_0) + \frac{C}{2}(x_1-x_0), \\
h(x_1) - \frac{C}{2}(x_1-x_0) &\leq& \frac{p(x_0)-p(x_1)}{x_0-x_1} &\leq& h(x_1) + \frac{C}{2}(x_1-x_0).
\end{aligned}
\]
Combining, we finally obtain
\[ -C \leq \frac{h(x_1)-h(x_0)}{x_1-x_0} \leq C.
\]
Hence $h(x)$ is locally Lipschitz and, in particular, continuous. Therefore $p \in C^1$. It also follows that $p_{\pm}(x)$ and its derivatives continuously depend on the initial point with respect to which these functions are defined.

For $x_1 > x_0$, $x_0,x_1 \in U$ we get
\begin{align}
    &p(x_0) + h(x_0)(x_1-x_0) + \frac{p_-''(\psi)}{2}(x_1-x_0)^2 = p_-(x_1) \leq p(x_1) \leq p_+(x_1) =\\\nonumber
    &p(x_0) + h(x_0)(x_1-x_0) + \frac{p_+''(\xi)}{2}(x_1-x_0)^2
\end{align}
where $\psi,\xi \in (x_0,x_1)$ are some intermediate points, or equivalently
\[ h(x_0) + \frac{p_-''(\psi)}{2}(x_1-x_0) \leq \frac{p(x_1)-p(x_0)}{x_1-x_0} \leq h(x_0) + \frac{p_+''(\xi)}{2}(x_1-x_0).
\]
By exchanging the role of $x_0,x_1$ we similarly get the inequalities
\[ h(x_1) - \frac{\tilde p_+''(\tilde\psi)}{2}(x_1-x_0) \leq \frac{p(x_0)-p(x_1)}{x_0-x_1} \leq h(x_1) - \frac{\tilde p_-''(\tilde\xi)}{2}(x_1-x_0),
\]
where $\tilde p_{\pm}$ now denote the corresponding functions with initial values $\tilde p_{\pm}(x_1) = p(x_1)$, $\tilde p'_{\pm}(x_1) = h(x_1)$, and $\tilde\psi,\tilde\xi \in (x_0,x_1)$ are other intermediate points. Combining, we obtain
\[ \frac{p_-''(\psi)+\tilde p_-''(\tilde\xi)}{2} \leq \frac{h(x_1)-h(x_0)}{x_1-x_0} \leq \frac{p_+''(\xi)+\tilde p_+''(\tilde\psi)}{2}.
\]
However, for $x_1 \to x_0 = x^*$ the left-most expression tends to 
\[ p_-''(x^*) = 6h(x^*)p(x^*) - 4p(x^*)^3 - 2\gamma(h(x^*) - p(x^*)^2)^{3/2},
\]
while the right-most expression tends to 
\[ p_+''(x^*) = 6h(x^*)p(x^*) - 4p(x^*)^3 + 2\gamma(h(x^*) - p(x^*)^2)^{3/2}.
\]
Similarly, for $x_0 \to x_1 = x^*$ the left-most expression tends to 
\[ \tilde p_-''(x^*) = 6h(x^*)p(x^*) - 4p(x^*)^3 - 2\gamma(h(x^*) - p(x^*)^2)^{3/2},
\]
while the right-most expression tends to 
\[ \tilde p_+''(x^*) = 6h(x^*)p(x^*) - 4p(x^*)^3 + 2\gamma(h(x^*) - p(x^*)^2)^{3/2}.
\]
This proves our claim.
\end{proof}

{\lemma \label{lem:sufficiency_boundary} Let $p,h: (-1,1) \to \mathbb R$ be functions as in Lemma \ref{lem:sufficiency2dim}. Then the function $f(x) = \int_0^x p(s)\,ds$ is admissible in the sense of Definition \ref{def:lipschitz}. It grows asymptotically proportional to $-\log(1-|x|)$ as $x \to \pm1$. }

\begin{proof}
By Lemma \ref{lem:sufficiency2dim} the function $f$ is $C^2$, defined on the whole interval $(-1,1)$, and satisfies the first and third property in Definition \ref{def:lipschitz}. It remains to show the second property.

From the conditions in Lemma \ref{lem:sufficiency2dim} it follows that
\[
\begin{aligned}
\frac{1}{\sqrt{\nu-1}}\left( -p(x) + \frac{1}{1-x} \right) &\leq & \sqrt{h(x)-p(x)^2} &\leq& \sqrt{\nu-1}\left( p(x) + \frac{1}{1+x} \right), \\
\frac{1}{\sqrt{\nu-1}}\left( p(x) + \frac{1}{1+x} \right) &\leq& \sqrt{h(x)-p(x)^2} &\leq& \sqrt{\nu-1}\left( -p(x) + \frac{1}{1-x} \right).
\end{aligned}
\]
Comparing the right-most with the left-most expressions and resolving with respect to $p(x)$ yields
\[ \frac{1}{\nu}\frac{1}{1-x} - \frac{\nu-1}{\nu}\frac{1}{1+x} \leq p(x) \leq \frac{\nu-1}{\nu}\frac{1}{1-x} - \frac{1}{\nu}\frac{1}{1+x}.
\]
Integrating we obtain
\[ -\frac{\nu-1}{\nu}\log(1+|x|)-\frac{1}{\nu}\log(1-|x|) \leq f(x) \leq -\frac{\nu-1}{\nu}\log(1-|x|)-\frac{1}{\nu}\log(1+|x|).
\]
The lower bound tends to $+\infty$ as $x \to \pm1$. This proves the second property in Definition \ref{def:lipschitz} and the claimed growth rate and completes the proof of the lemma.
\end{proof}

Finally, let us eliminate the dependence on the variable $x \in (-1,1)$ in the conditions of Lemma \ref{lem:sufficiency2dim}. We consider the slightly more restrictive, but simpler case of $C^3$ functions, as in the original definition of self-concordance. 

\begin{lemma} \label{lem:centeredCondition}
Let $f: (-1,1) \to \mathbb R$ be a $C^3$ function, and let $x \in (-1,1)$ be arbitrary. Define $p = f'(x)$, $h = f''(x)$, $w = f'''(x)$ and
\[
\tilde p = (1-x^2)p - x,\quad \tilde h = (1-x^2)^2h - 2x(1-x^2)p + x^2,
\]
\[ \tilde w = (1-x^2)^3w - 6x(1-x^2)^2h + 6x^2(1-x^2)p - 2x^3.
\]
Then conditions \eqref{lemma4_1},\eqref{lemma4_2} are equivalent to the conditions
\begin{equation} \label{Cnu2condition}
    \frac{\sqrt{\tilde h-\tilde p^2}}{\sqrt{\nu-1}} \leq \pm\tilde p + 1 \leq \sqrt{\nu-1}\sqrt{\tilde h-\tilde p^2},
\end{equation}
and the condition
\[
|w - 6hp + 4p^3| \leq 2\gamma(h - p^2)^{3/2}
\]
is equivalent to the condition
\begin{equation} \label{Cnu3condition}
    |\tilde w - 6\tilde h\tilde p + 4\tilde p^3| \leq 2\gamma(\tilde h - \tilde p^2)^{3/2}.
\end{equation}
\end{lemma}

\begin{proof}
    The proof is by direct substitution, passing by the intermediate steps
    \[
    \tilde h - \tilde p^2 = (1-x^2)^2(h - p^2), \quad \tilde w - 6\tilde h\tilde p + 4\tilde p^3 = (1-x^2)^3(w - 6hp + 4p^3).
    \]
\end{proof}

We obtain the following result.

\begin{corollary} \label{1DequivalenceBody}
    Let $f: (-1,1) \to \mathbb R$ be a $C^3$ function. For every $x \in (-1,1)$, define $\tilde p,\tilde h,\tilde w$ as in Lemma \ref{lem:centeredCondition}. Then the following conditions are equivalent:
    \begin{itemize}
        \item the function $f$ satisfies the conditions in Lemma \ref{lem:section}
        \item for every $x$ the triple $(\tilde p,\tilde h,\tilde w)$ is an element of the set $P_{\nu}$ defined by conditions \eqref{Cnu2condition},\eqref{Cnu3condition}, i.e.,
        \[ P_{\nu} = \left\{ (x_1,x_2,x_3) \mid \frac{\sqrt{x_2-x_1^2}}{\sqrt{\nu-1}} \leq \pm x_1 + 1 \leq \sqrt{\nu-1}\sqrt{x_2-x_1^2},\right.
        \]
        \[ \left. |x_3 - 6x_2x_1 + 4x_1^3| \leq 2\frac{\nu-2}{\sqrt{\nu-1}}(x_2-x_1^2)^{3/2} \right\}.
        \]
    \end{itemize}
\end{corollary}

The proof is by combination of the results above.

In this section we obtained a characterization of the admissibility of the function $f: (-1,1) \to \mathbb R$ in the sense of Lemma \ref{lem:centeredCondition} by inequalities on the values of the derivatives $f'(x),f''(x),f'''(x)$. These conditions consist in the inclusion of an affine image of the triple of derivatives into a certain non-convex body $P_{\nu}$. It is important to note that the affine image depends on the point $x \in (-1,1)$ but the body $P_{\nu}$ does not.

\subsection{Cones of arbitrary dimension}

In this section we generalize the results of the previous section to the case of compact sets $C$ of arbitrary finite dimension.

{\lemma \label{lem:gendim_admissible} Let $C \subset \mathbb R^n$ be a compact convex set. A $C^2$ function $f: C^o \to \mathbb R$ is admissible in the sense of Definition \ref{def:lipschitz} if and only if its restrictions to every interval $I \subset C^o$ linking two boundary points of $C$ are admissible. }

\begin{proof}
Clearly if $f$ is admissible, then so are its restrictions to interior intervals linking boundary points of $C$.

Let us prove the reverse implication. The first and third property in Definition \ref{def:lipschitz} for $f$ on $C^o$ follow immediately from the corresponding properties for the restrictions of $f$ to intervals. The second property follows from the logarithmic growth rate of $f$ along any straight line when approaching a boundary point, which was established in Lemma \ref{lem:sufficiency_boundary}.
\end{proof}

We conclude the following result.

\begin{theorem} \label{thm:bodiesCharacterization}
Let $C$ be a compact convex set and let $f: C^o \to \mathbb R$ be a $C^3$ function on the interior of $C$. For every interior point $x \in C^o$ and every non-zero direction $v$, let $I$ be the interval obtained by the intersection of $C$ with the line $x + \mathbb R \cdot v$. Parameterize this interval affinely by $t \in [-1,1]$ and consider the restriction $f_I: (-1,1) \to \mathbb R$ of $f$ on $I$. Define $p = f_I'(x)$, $h = f_I''(x)$, $w = f_I'''(x)$ to be the first three derivatives with respect to the parameter $t$ at the point $x$. Let the triple $(\tilde p,\tilde h,\tilde w)$ be defined as in Lemma \ref{lem:centeredCondition}. 

Then the following conditions are equivalent:
\begin{itemize}
  \item for all $x \in C^o$ and all non-zero directions $v$, the triple $(\tilde p,\tilde h,\tilde w)$ is an element of the set $P_{\nu}$
  \item the function $f$ satisfies the conditions of Lemma \ref{lem:section} on $C$.
\end{itemize}
\end{theorem}

\begin{proof}
    Suppose that the triples $(\tilde p,\tilde h,\tilde w)$ are in $P_{\nu}$ for all interior points $x$ and all directions $v$. By Corollary \ref{1DequivalenceBody} the restriction $f_I$ satisfies the conditions of Lemma \ref{lem:section} on the interior of every interval $I$ linking two boundary points of $C$. By Lemma \ref{lem:gendim_admissible} the function $f$ then also satisfies the conditions of Lemma \ref{lem:section}.

    The conclusion in the reverse direction follows similarly.
\end{proof}

Theorem \ref{thm:bodiesCharacterization} has the following consequence.

\begin{corollary}
    Let $f_1,\dots,f_k$ be $C^3$ functions satisfying the conditions of Lemma \ref{lem:section} with parameter $\nu$ on some convex compact set $C$. Let $f_c$ be a convex combination of $f_1,\dots,f_k$.

    Suppose that $\tilde{\nu}$ is such that $P_{\tilde\nu}$ contains the convex hull of the body $P_{\nu}$, where the family of bodies $P_{\nu}$ is defined in Corollary \ref{1DequivalenceBody}. Then $f_c$ satisfies the conditions of Lemma \ref{lem:section} with parameter $\tilde{\nu}$.
\end{corollary}

\begin{proof}
    By Theorem \ref{thm:bodiesCharacterization} the triples $(\tilde p_i,\tilde h_i,\tilde w_i)$ corresponding to the function $f_i$ are in $P_{\nu}$ for all interior points $x \in C^o$ and all non-zero directions $v$.

    However, the triple $(\tilde p,\tilde h,\tilde w)$ depends affinely on the function $f$. Hence the triple $(\tilde p_c,\tilde h_c,\tilde w_c)$ corresponding to the convex combination $f_c$ is in the convex hull of $P_{\nu}$ for all interior points $x$ and all directions $v$. But then $(\tilde p_c,\tilde h_c,\tilde w_c) \in P_{\tilde{\nu}}$, and it follows that $f_c$ satisfies the conditions of Lemma \ref{lem:section} with parameter $\nu$.
\end{proof}

Note that the gap between $\nu$ and $\tilde\nu$ is small for values of $\nu$ close to 2. Since the results of this paper are anyway intended for computing barriers on 3- or 4-dimensional cones due to the otherwise prohibitive computational complexity, and hence for small values of the optimal parameter $\nu$, this gap is only a minor drawback.

The goal of this paper is to quantify this degradation of the parameter from the value $\nu$ to the value $\tilde{\nu}$. To this end we have to establish, for given $\nu$, the least $\tilde{\nu}$ such that $P_{\tilde{\nu}}$ contains the convex hull of $P_{\nu}$.

\section{Numerical calculation of $\tilde{\nu}$}\label{sec:numerical}

\begin{figure}[h]
\centering
 \begin{subfigure}[c]{0.15\textwidth}
    \caption{$\nu = 2.5$}  \includegraphics[width=5.00cm,height=2.50cm]{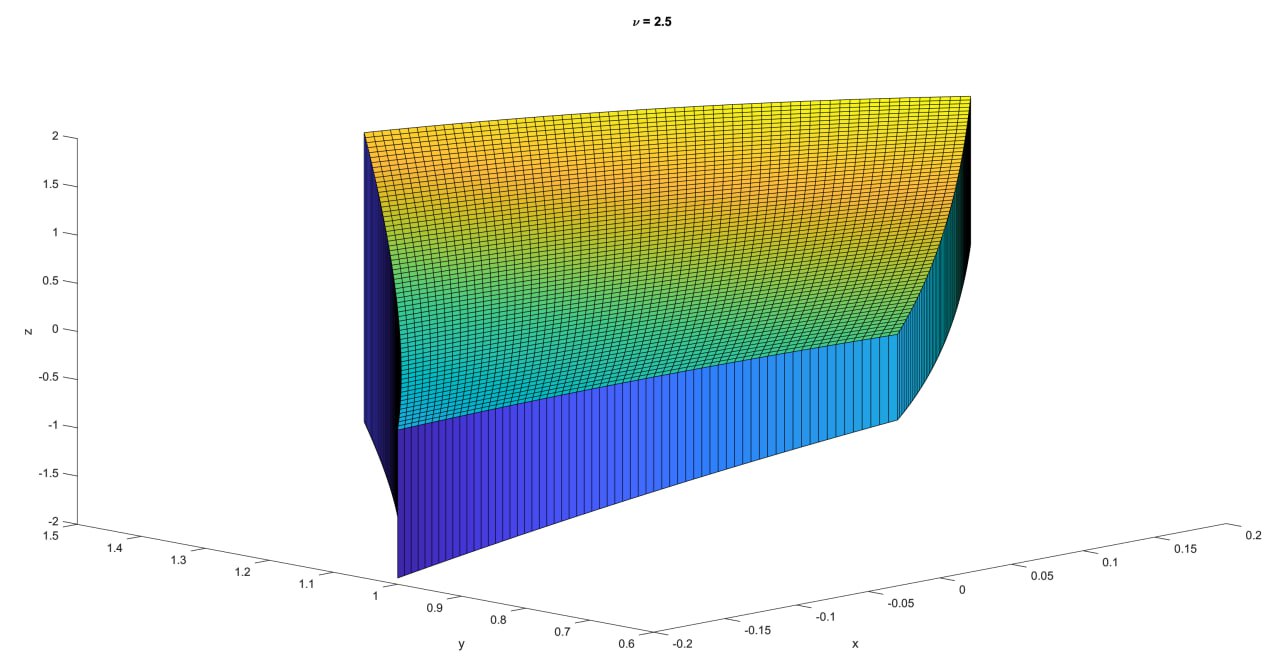}
\end{subfigure}\hfill
 \begin{subfigure}[c]{0.45\textwidth}
    \caption{$\nu = 3$}  \includegraphics[width=5.00cm,height=2.50cm]{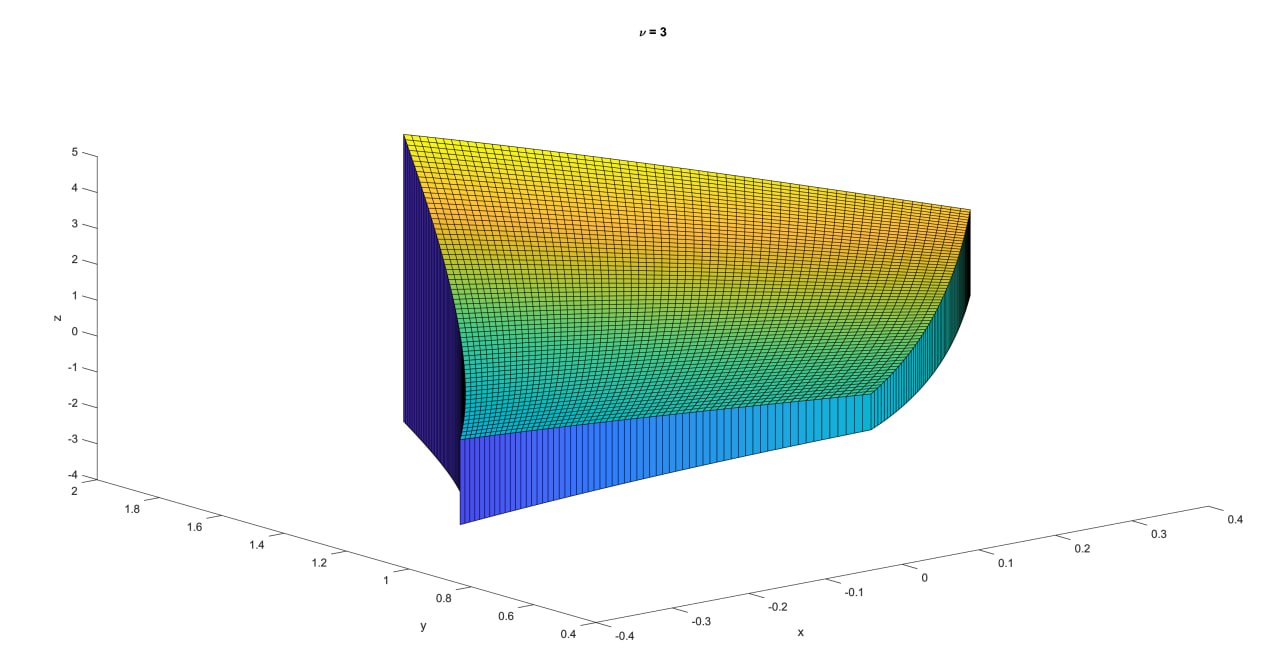}
\end{subfigure}\hfill
 \begin{subfigure}[c]{0.45\textwidth}
    \caption{$\nu = 4$}  \includegraphics[width=5.00cm,height=2.50cm]{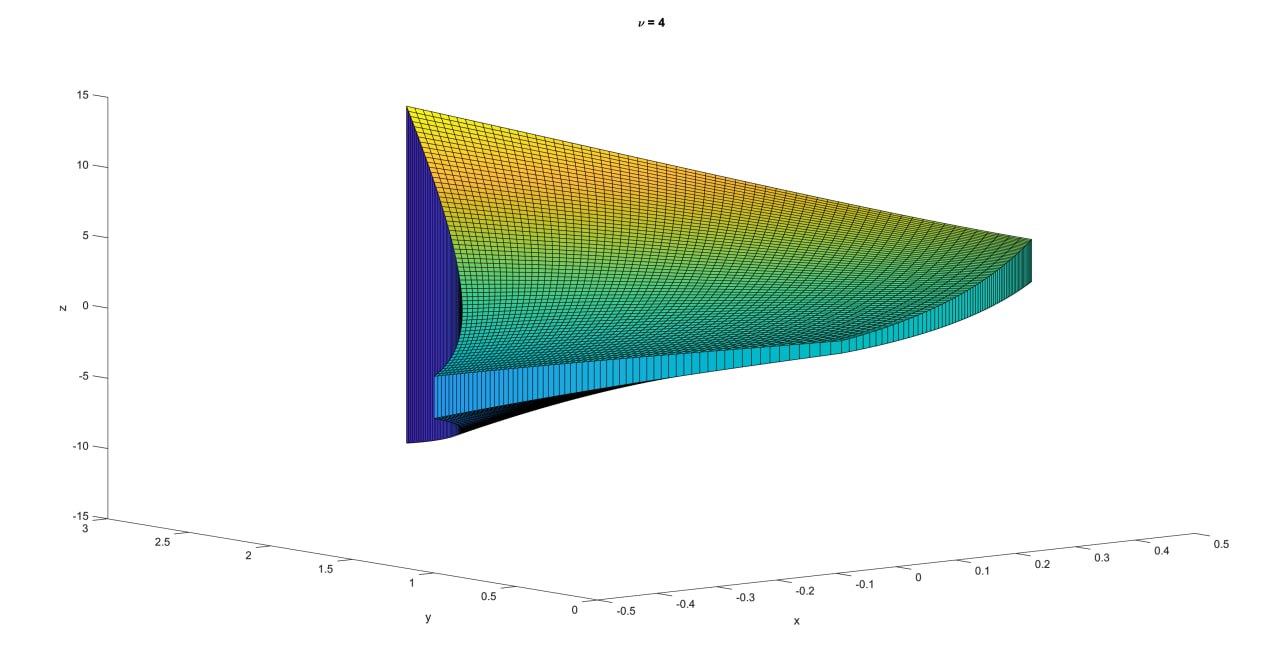}
\end{subfigure}\hfill

\caption{$P(\nu)$ set for different $\nu$ (2.5, 3, and 4 respectively)}
\label{fig:body}
\end{figure}

We have numerically calculated the body $P(\nu)$ and depicted it on Fig.~\ref{fig:body} for several values of $\nu$. The calculation was conducted on a uniform grid with 101 nodes, with $(x_1, x_2)$ calculated on a grid with a step size of approximately 0.012. The calculation was conducted for $\nu$ in the range $[2.1, 5.1]$. It was found that the set $P(\nu)$ becomes more curved with increasing value of the parameter, resulting in a less pronounced resemblance to its convex hull. This phenomenon leads to a fast increase in the function $\tilde{\nu}$ for large $\nu$.

\begin{figure}[h]
\centering
 \includegraphics[width=7.00cm,height=5.00cm]{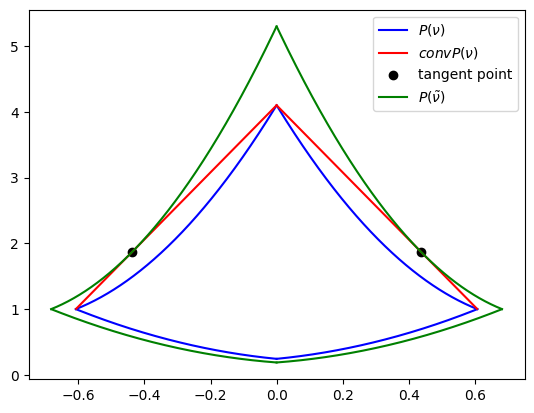}
\caption{Projections on the $x_1,x_2$ plane of the set $P(\nu)$ for $\nu = 5.1$ and $\tilde{\nu} = 6.30125$.}
\label{fig:projection}
\end{figure}

A lower bound on $\tilde{\nu}$ can be calculated by determining the projection of the convex hull of $P(\nu)$ onto the $(x_1, x_2)$ space, and subsequently finding the lowest possible $\tilde{\nu}$ for which the projection is contained in $P(\tilde{\nu})$. The projection is symmetric with respect to the reflection $x_1 \mapsto -x_1$ and bounded by four parabolas, two of which are convex and two concave. The upper right concave boundary arc of the projection is given by the graph of the function $x_2 = (\nu - 1)(x_1 - 1)^2 + x_1^2$ on the interval $x_1 \in [0,\frac{\nu-2}{\nu}]$. A direct calculation shows that its convex hull is bounded by the line $x_2 = -\nu x_1 + \nu - 1$. This line is tangent to the graph of the function $x_2 = (\nu' - 1)(x_1 - 1)^2 + x_1^2$ with $\nu' = \frac{\nu^2}{8} + \frac{\nu}{2} + \frac{1}{2}$. It can thus be concluded that $\tilde{\nu} \geq \frac{\nu^2}{8} + \frac{\nu}{2} + \frac{1}{2}$. The projection for $\nu = 5.1$ and $\nu = 6.30125$ is presented in Fig.~\ref{fig:projection}.

Using the moment techniques presented in \cite{lasserre2009moments}, we have constructed the convex hull $\conv\,P(\nu)$ numerically on the same grid, using semi-definite programming to describe the convex hull of rational curve segments. Subsequently, for each value of the parameter $\nu$ and for each point $(x_1, x_2)$ in the grid, the minimum possible value of the parameter $\tilde{\nu}$ for which the upper and lower limits of the $x_3$ coordinate of the point in question lie between the extreme values of the $x_3$ coordinate in the set $P(\tilde{\nu})$ are determined. The maximum at each point is assumed to be a good estimate of the true maximum.

\begin{figure}[h]
\centering
 \includegraphics[width=12.00cm,height=5.00cm]{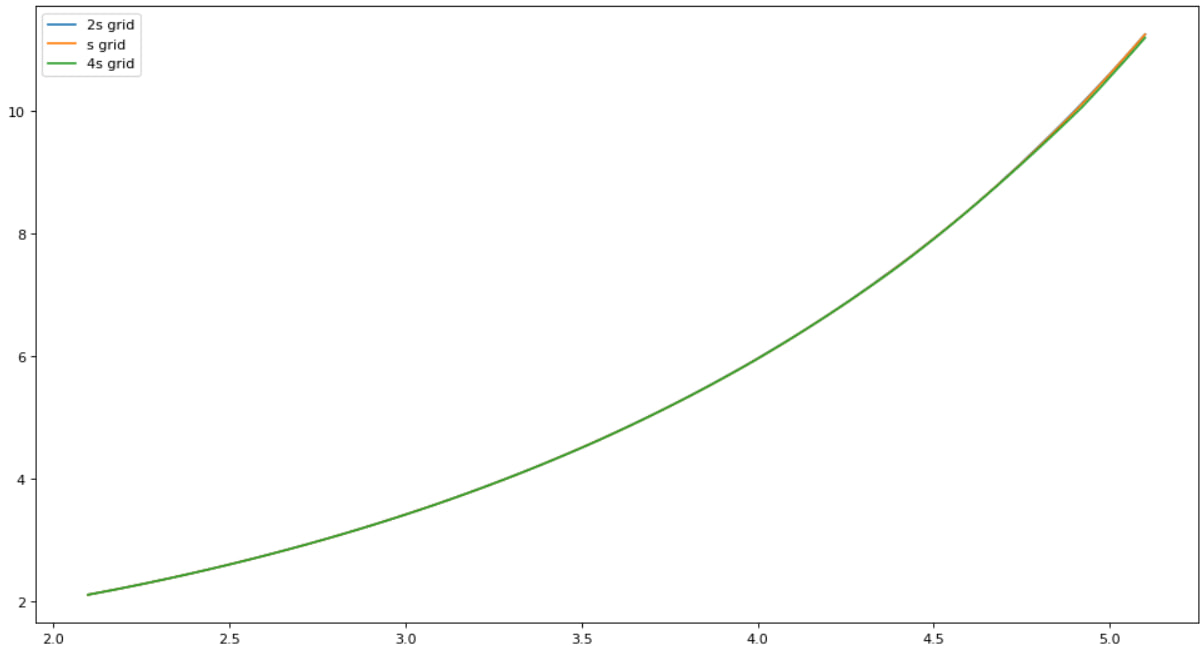}
\caption{Value $\tilde{\nu}$ calculated for 101 $\nu$ values uniformly distributed between 2.1 and 5.1 on different ($x_1, x_2$) grids}
\label{fig:tilde_nu}
\end{figure}

\begin{figure}[h]
\centering
 \includegraphics[width=12.00cm,height=5.00cm]{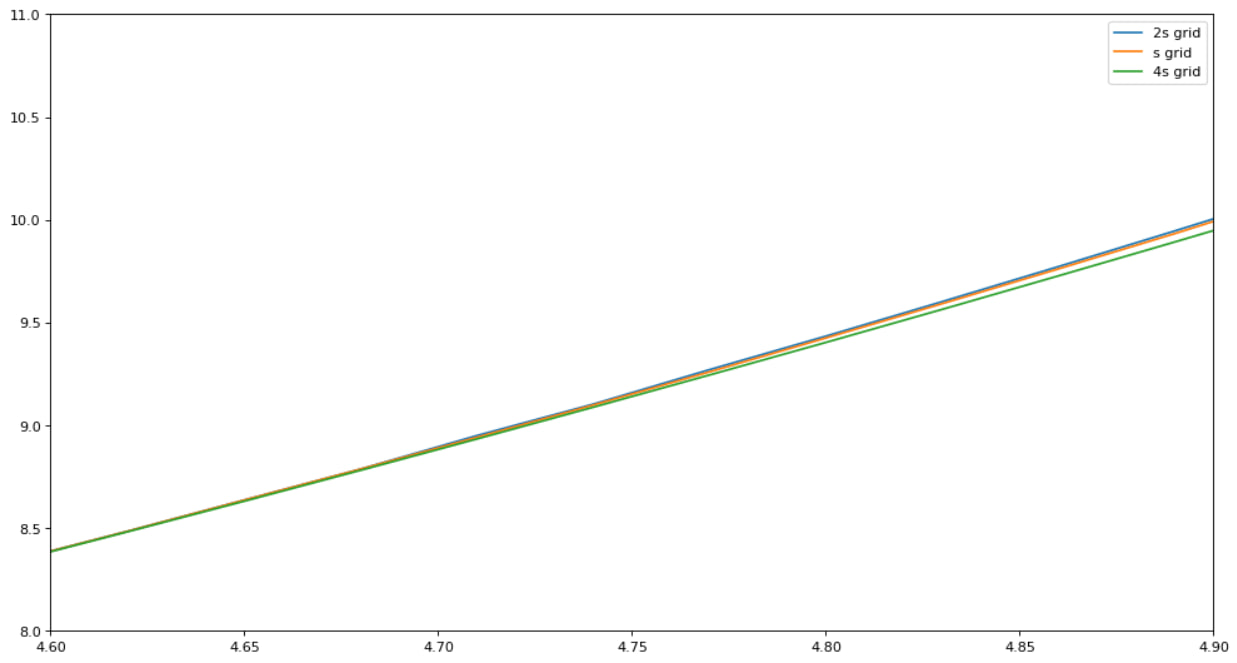}
\caption{Value $\tilde{\nu}$ calculated for 101 $\nu$ values uniformly distributed between 2.1 and 5.1 on different ($x_1, x_2$) grids}
\label{fig:tilde_nu_zoom}
\end{figure}

The accuracy of this estimate has been determined by means of a straightforward method. It can be shown that, due to the continuity of the set family, the true value of the bound $\tilde\nu$ is equal to the limit of the calculated values on the grid as the step approaches zero. To this end, the value of the function for three distinct grids with steps $s = 0.012$, $2s$, and $4s$ has been calculated, and the degradation has been examined. Due to continuity, the degradation of a grid with step $s$ compared to a perfect grid is comparable to the degradation of a grid with step $4s$ or $2s$ compared to a grid with step $s$. The results are presented in figure \ref{fig:tilde_nu}. To illustrate the difference more clearly, we have included a zoomed version of the plot (see fig. \ref{fig:tilde_nu_zoom}). The absolute values of the difference between the $4s$ grid and the $s$ grid is 0.5, and between the $2s$ grid and the $s$ grid is 0.09.


\section{Conclusion}
In this paper we have formulated necessary conditions for the existence of a barrier with given parameter $\nu$ and have constructed their convexification. Unfortunately, the convexification leads to a parameter degradation, but this degradation can be quantified. More precisely, a value $\nu$ corresponding to the convex relaxation degrades to a true value $\nu + O(\nu - 2)^2$. We leave the question of convexifying sufficient conditions open for further research.

\bibliographystyle{plain}
\bibliography{references}

\end{document}